\documentclass[12pt]{article}
\usepackage{amssymb,amsfonts,amsthm}
\usepackage{eucal,latexsym}
\newcounter{ppp}

\newcommand{\la}{\langle}
\newcommand{\ra}{\rangle}

\newcommand{\mod}{{\hbox{mod }}}

\newcommand{\ttt}{{\cal T}}

\newcommand{\Ker}{\hbox{Ker}}

\newcommand{\iv}{^{-1}}

\newcommand{\F}{\mathbb{F}}
\newcommand{\Fr}{\hbox{Fr}}
\newcommand{\HNN}{\hbox{HNN}}

\newcommand{\Z}{{\mathbb Z}}

\newcommand{\SL}{\hbox{SL}}
\newcommand{\GL}{\hbox{GL}}
\newcommand{\PGL}{\hbox{PGL}}
\newcommand{\Spec}{{\rm Spec}}

\newcommand{\adj}{{\hbox{adj}}}
\begin{document}
\renewcommand{\theequation}{\thesection.\arabic{equation}}
\newtheorem{theo}{\quad Theorem}[section]
\newtheorem{lemma}[theo]{\quad Lemma}
\newtheorem{cy}[theo]{\quad Corollary}
\newtheorem{df}[theo]{\quad Definition}
\newtheorem{rk}[theo]{\quad Remark}
\newtheorem{prop}[theo]{\quad Proposition}
\newtheorem{prob}[theo]{\quad Problem}
\newtheorem{conj}[theo]{\quad Conjecture}

\title{Polynomial maps over finite fields and
residual finiteness of mapping tori of group endomorphisms}
\author{Alexander Borisov and Mark Sapir\thanks{The research of the
second author was supported in part by the NSF grants DMS 9978802,
0072307, 0245600, and the US-Israeli BSF grant 1999298.}}
\date{}
\maketitle

\begin{abstract} We prove that every mapping torus of any
free group endomorphism is residually finite. We show how to use a
not yet published result of E. Hrushovski to extend our result to
arbitrary linear groups. The proof uses algebraic self-maps of
affine spaces over finite fields. In particular, we prove that
when such a map is dominant, the set of its fixed closed scheme
points is Zariski dense in the affine space.
\end{abstract}

\section{Introduction}

This article contains results in group theory and algebraic
geometry. We think that both the results and the relationship
between them are interesting and will have other applications in
the future.

We start with group theory. Let $G$ be a group given by generators
$x_1,...,x_k$ and a set of defining relations $R$, and let
$\phi\colon x_i\mapsto w_i$, $1\le i\le k$ be an injective
endomorphism of $G$. Then the group $$\HNN_\phi(G)=\la
x_1,...,x_k, t\mid R, tx_it\iv=w_i, i=1,...,k\ra$$ is called the
{\em mapping torus} of $\phi$ (or {\em ascending HNN extension of
$G$ corresponding to} $\phi$). This group has an easy geometric
interpretation as the fundamental group of the mapping torus of
the standard 2-complex of $G$ with bounding maps the identity and
$\phi$. The simplest and one of the most important cases is when
$G$ is the free group $F_k$ of rank $k$, i.e. when $R$ is empty.
These groups appear often in group theory and topology. In
particular, many one-relator groups are ascending HNN extensions
of free groups (more on that below).

Some essential information about the mapping tori of free group
endomorphisms is known. In particular, Feighn and Handel \cite{FH}
proved that these groups are coherent, that is, all their finitely
generated subgroups are finitely presented. They also
characterized all finitely generated subgroups of such groups. We
know \cite{GMSW} that these groups are Hopfian, that is every
surjective endomorphism of such a group must be injective. On the
other hand, ascending HNN extensions of arbitrary residually
finite groups are not necessarily Hopfian \cite{SW}.

Many of the groups of the form $\HNN_\phi(F_k)$ are hyperbolic
(see \cite{BF} and \cite{Kapovich}). One of the outstanding
problems about hyperbolic groups is whether they are residually
finite. Recall that a group is called {\em residually finite} if
the intersection of its subgroups of finite index is trivial. This
leads to the following question:

\begin{prob} \label{pr1} Are all mapping tori of free groups residually finite?
\end{prob}

This question also arises naturally when one tries to characterize
residually finite one-related groups. As far as we know Problem
\ref{pr1} was explicitly formulated first by Moldavanskii in
\cite{Mol} (it is also mentioned in \cite{Wise1} and listed as
Problem 1 in the list of ten interesting open problems about
ascending HNN extensions of free groups in \cite{Kapovich}).

Notice that ascending HNN extensions of residually finite groups
may be not residually finite. They can even have very few finite
homomorphic images as is the case for Grigorchuk's group
\cite{SW}. However if $\phi$ is an automorphism and $G$ is
residually finite then $\HNN_\phi(G)$ is also residually finite
\cite{Mal1}. Thus the interesting case in Problem \ref{pr1} is
when $\phi$ is not surjective. Some special cases of Problem
\ref{pr1} have been solved in \cite{Wise1} (these cases proved to
be useful in Wise's residually finite version of Rips'
construction), and in \cite{HW} (where it is proved that the
mapping tori of polycyclic groups are residually finite).

Notice also that the groups $\HNN_\phi(F_k)$ do not necessarily
satisfy properties that are known to be somewhat stronger than the
residual finiteness. For example, the group $\la a, t\mid
tat\iv=a^2\ra$ is not a LERF group ($a$ cannot be separated from
the cyclic group $\la a^2\ra$ by a homomorphism onto a finite
group).

One of the main goals of this paper is to solve Problem \ref{pr1}.

\begin{theo} \label{th} The mapping torus of any injective endomorphism of
a free group is residually finite.
\end{theo}

Computer experiments conducted by Ilya Kapovich, Paul Schupp, and
the second author of this paper seem to show that most 1-related
groups are subgroups of ascending HNN extensions of a free
group\footnote{A simple Maple program written by the second author
of this paper checked 30,000 random two-letter group words of
length 300,000 Schupp's program checked 50,000 two-letter random
words of length between 100,000 and 110,000. Both programs found
that at least 99.6\% of the corresponding 1-related groups are
subgroups of ascending HNN extensions of finitely generated free
groups.}. Thus it could well be true that groups with one defining
relation are generically inside ascending HNN extensions of free
groups. If this conjecture turns out to be true then Theorem
\ref{th} would imply that one-related groups are generically
residually finite. (Recall that there exist non-residually finite
one-related groups, for example the Baumslag-Solitar group
$BS(2,3)=\la a,t\mid ta^2t\iv=a^3\ra$.) Anyway, it is clear that
Theorem \ref{th} applies to very many one-related groups.

The proof of Theorem \ref{th} was obtained in a rather unexpected
way. The proofs of the previous major results about mapping tori
of groups (see for example \cite{FH}, \cite{GMSW},
\cite{Kapovich}) were of topological nature. We know of several
attempts (see \cite{HW}, \cite{Wise1}) to apply similar methods to
Problem \ref{pr1}: residual finiteness of the fundamental group of
a CW-complex is equivalent to the existence of enough finite
covers of that complex to separate all elements of the fundamental
groups. But these approaches produced only partial results. Even
simple examples like the group $\la a, b, t\mid tat\iv=ab,
tbt\iv=ba\ra$ have been untreatable so far by the topological
methods.

Our approach is based on a reduction of Problem \ref{pr1} to some
questions about periodic orbits of algebraic maps over finite
fields (see Section 2). More precisely, we study the orbits
consisting of points conjugate over the base field. In the
language of schemes these orbits correspond to the fixed closed
scheme points. Such points appeared in the Deligne Conjecture, and
were extensively studied before (see, e.g., \cite{Fu},
\cite{Pink}). However, these investigations were limited to the
quasi-finite maps and  most of our maps are not quasi-finite.

\begin{df} Suppose $\Phi\colon X\to X$ is a self-map of a variety over a finite
field $\F_q$. A geometric point $x$ of $X$ over some finite
extension of $\F_q$ is called quasi-fixed with respect to $\Phi$
if $\Phi(x)=\Fr^m(x)$. Here $\Fr^m$ is the m-th composition power
of the geometric Frobenius morphism.
\end{df}

If $X$ is the affine space, the above definition becomes the following.
Let $\Phi\colon A^n(\F_q)\to A^n(\F_q)$ be a polynomial map,
defined over the finite field $\F_q.$ It is given in coordinates
by the polynomials $\phi_1,...,\phi_n$ from $\F_q[x_1,...,x_n]$. Suppose
a point $a=(a_1,a_2,...,a_n) \in A^n$ is defined over the
algebraic closure $\bar{\F}_q$ of $\F_q.$ It is a
quasi-fixed point of $\Phi$ if and only if for some $Q=q^m$ for all $i$
$$\phi_i(a_1,a_2,...,a_n)=a_i^Q.$$

Here is our main theorem regarding such maps.

\begin{theo}\label{th0}
Let $\Phi^{n}\colon A^n(\F_q)\to A^n(\F_q)$ be the $n$-th iteration
of $\Phi$. Let $V$ be the Zariski closure of $\Phi^{n}(A^n)$. Then
the following holds.
\begin{enumerate}
\item All quasi-fixed points of $\Phi$ belong to $V$.

\item Quasi-fixed points of $\Phi$ are Zariski dense in $V$. In
other words, suppose $W\subset V$ is a proper Zariski closed
subvariety of $V$. Then for some $Q=q^m$ there is a point
$(a_1,...,a_n)\in U\setminus W$ such that for all $i$
$f_i(a_1,...a_n)=a_i^Q.$
\end{enumerate}
\end{theo}

After we obtained the proof of Theorem \ref{th0}, we received a
preprint \cite{Hr} of E. Hrushovski where he proves a more general
result. In particular, his results imply the following

\begin{theo}\label{Hr}(Hrushovski, \cite{Hr})
Let $\Phi\colon X\to X$ be a dominant self-map of an absolutely irreducible variety over a finite field. Then the set of the quasi-fixed points of $\Phi$ is Zariski dense in $X.$
\end{theo}

Our Theorem \ref{th0} is a partial case of Theorem \ref{Hr} where
$X$ is the Zariski closure of $\Phi^n(A^n)$ (in particular, our
theorem captures the case when $\Phi\colon A^n\to A^n$ is
dominant).

Theorem \ref{Hr} allowed us to prove the following statement that
is much stronger than Theorem \ref{th}.

\begin{theo}\label{thH} The mapping torus of any injective endomorphism
of a finitely generated linear group\footnote{That is a group
representable by matrices of any size over any field.} is
residually finite.
\end{theo}

As we mentioned before, for non-linear residually finite groups
this statement is not true \cite{SW}. In fact Theorem \ref{th} can
serve as a tool to show that a group is {\em not} linear. For
example, the non-Hopfian example from \cite{SW} is an ascending
HNN extension of a residually finite finitely generated group that
is an amalgam of two free groups. By Theorem \ref{thH} that
amalgam of free groups is not linear.

It is well known that free groups, polycyclic groups, etc. are
linear. Thus Theorem \ref{th} immediately implies all known
positive results about residual finiteness of mapping tori of
non-surjective endomorphisms \cite{HW}, \cite{GMSW}, \cite{Mol},
\cite{Wise1}.

The proof of Theorem \ref{Hr} is complicated and uses some heavy
machinery from algebraic geometry and Hrushovski's theory of
difference schemes. In comparison, our proof of Theorem \ref{th0}
is basically elementary.

\begin{rk} {\rm Theorems \ref{th} and \ref{thH} will remain true if we
drop the requirement that the endomorphism $\phi$ is injective.
Indeed, it is easy to see that for every endomorphism $\phi$ of a
linear group $G$, the sequence $\Ker(\phi)\subseteq
\Ker(\phi^2)\subseteq \Ker(\phi^3)\subseteq...$ eventually
stabilizes (see \cite[Theorem 11]{Mal}). Then, for some $n$,
$\phi$ is injective on $\phi^n(G)$, and the group $HNN_\phi(G)$ is
isomorphic to the ascending HNN extension of $\phi^n(G)$. Since
$\phi^n(G)$ is again a linear group, we can apply Theorem
\ref{thH} (see details in \cite{Kap02}).}
\end{rk}

The paper is organized as follows. In Section 2, we reduce Theorems
\ref{th} and \ref{thH} to Theorems \ref{th0} and \ref{Hr}. In
Section 3, we give a proof of Theorem \ref{th0}. In Section 4 we
apply Theorem \ref{thH} to a question about extendability of
endomorphisms of linear groups to automorphisms of their profinite
completions. We also present some open problems.

{\bf Acknowledgments.} The authors are grateful to Ilya Kapovich,
Yakov Varshavsky and Dani Wise for very useful remarks.

\section{HNN extensions and dynamical systems}

Let $T=\HNN_\phi(G)$ be the ascending HNN extension of a group
$$G=\la x_1,...,x_k\mid R\ra$$ corresponding to an injective
endomorphism $\phi$. Let $t$ be the free letter of this HNN
extension, so that $tx_it\iv = \phi(x_i)$ for every $i=1,...,k$.

It is easy to see that every element $g$ of $T$ can be written as
a product $t^awt^b$ for some integers $a\le 0$ and $b\ge 0$, $w\in
G$. The map $z: T\to \Z$ that sends $t^awt^b$ to $a+b$ is a
homomorphism, so if $a+b\ne 0$ then $g$ can be separated from 1 by
a homomorphism onto a finite group. If $a=-b$ then $g$ and $w$ are
conjugate, so for every homomorphism $\psi$, $\psi(g)\ne 1$ if and
only if $\psi(w)\ne 1$. Therefore the following fact is true.

\begin{lemma}\label{lm1}
The group $T$ is residually finite if and only if for every $w\in
G$, $w\ne 1$, there exists a homomorphism $\psi$ of $T$ onto a
finite group such that $\psi(w)\ne 1$.
\end{lemma}

Let $\phi$ be an endomorphism of $G$ defined by a sequence of
words $w_1,...,w_k$ from $F_k$ (that is the images of $w_i$ in $G$
under the natural homomorphism $F_k\to G$ generate a subgroup that
is isomorphic to $G$). Let $H$ be any group (or, more generally, a
group scheme). Then we can define a map $\phi_H\colon H^k\to H^k$
that takes every $k$-tuple $(h_1,...,h_k)$ to the $k$-tuple
$$(w_1(h_1,...,h_k), w_2(h_1,...,h_k), \ldots,w_k(h_1,...,h_k)).$$
Notice that this map is not a homomorphism. Nevertheless it
defines a dynamical system on $H^k$.

 The following lemma
reformulates residual finiteness in terms of these dynamical
systems.

\begin{lemma}\label{lm2} The group $T=\HNN_\phi(G)$ is residually finite if and
only if for every $w=w(x_1,...,x_k)\ne 1$ in $G$ there exists a
finite group $H=H_w$ and an element $h=(h_1,...,h_k)$ in $H^k$
such that
\begin{enumerate}

\item[(i)] $h_1,...,h_k$ satisfy the relations from $R$ (where
$h_i$ is substituted for $x_i$, $i=1,...,k$).

\item[(ii)] $h$ is a fixed point of some power of $\phi_H$, and

\item[(iii)] $w(h_1,...,h_k)\ne 1$ in $H$.
\end{enumerate}
\end{lemma}

\proof $\Rightarrow$ Suppose $T$ is residually finite. Take any
word $w\ne 1$ in $G$. Then there exists a homomorphism $\gamma$
from $G$ onto a finite group $H$ such that $\gamma(w)\ne 1$. Let
$t$ be the free letter in $G$. Then
$\gamma(t)\gamma(G)\gamma(t\iv)\subseteq \gamma(G)$. Since $H$ is
finite, $\gamma(t)$ acts on $\gamma(G)$ by conjugation. It is
clear that for every element $h=(h_1,...,h_k)$ in $\gamma(G)^k$,
\begin{equation}\phi_H(h)=(\gamma(t)h_1\gamma(t\iv),...,\gamma(t)h_k\gamma(t\iv))\in
\gamma(F_k)^k.\label{eq1}\end{equation} Take
$h=(\gamma(x_1),...,\gamma(x_k))$. Property (i) is obvious.
Property (iii) holds because $\gamma(w)\ne 1$. Property (ii) holds
also because by (\ref{eq1}) powers of $\phi_H$ act on
$\gamma(G)^k$ as conjugation by the corresponding powers of
$\gamma(t)$, and some power of $\gamma(t)$ is equal to 1 since $H$
is finite.

$\Leftarrow$ Suppose that for every $w\ne 1$ in $G$ there exists a
finite group $H=H_w$ and an element $h=(h_1,...,h_k)$ in $H^k$
such that conditions (i), (ii) and (iii) hold. We need to prove
that $G$ is residually finite. By Lemma \ref{lm1}, it is enough to
show that every such $w$ can be separated from 1 by
a homomorphism of $G$ onto a finite group.

Pick a $w\ne 1$ in $G$. Let a finite group $H$,  $h\in H^k$, be as
above. By (ii), there exists an integer $n\ge 1$ such that
$\phi_H^n(h)=h$. Let $P$ be the wreath product of $H$ and a cyclic
group $C=\la c\ra$ of order $n$. Recall that $P$ is the semidirect
product of $H^n$ and $C$ where elements of $C$ act on $H^n$ by
cyclically permuting the coordinates.

Consider the $\phi_H$-orbit $h^{(0)}=h, h^{(1)}=\phi_H(h), ...,
h^{(n-1)}=\phi_H^{n-1}(h)$ of $h$. Let
$h^{(i)}=(h_1^{(i)},...,h_k^{(i)})$, $i=0,...,n-1$. For every
$j=1,...,k$ let $y_j$ be the $n$-tuple
$(h_j^{(0)},h_j^{(1)},...,h_j^{(n-1)})$. Notice that since $h$
satisfies relations from $R$, $\phi_H(h)$, $\phi_H^2(h)$,... also
satisfy these relations because $\phi$ is an endomorphism of $G$.
This and (ii) immediately imply that the map $\phi\colon t\mapsto
c, x_j\mapsto y_j, j=1,...,k,$ can be extended to a homomorphism
of $T$ onto a subgroup of $P$ generated by $c, y_1,...,y_k$.
Notice that the image of $w$ under this homomorphism is an
$n$-tuple $w(y_1,...,y_k)$ from $H^n$ whose first coordinate is
$w(h_1,...,h_k)\ne 1$ in $H$ by property (iii). Thus $w$ can be
separated from 1 by a homomorphism of $T$ onto a finite group.
\endproof

Now we are going to show how to apply Lemma \ref{lm2} to free
groups and other linear groups. First we need to fix some
notations.

Let us identify the scheme $M_r$ of all $r$ by $r$ matrices with
the scheme $\Spec\Z[a_{i,j}]$, $1\le i,j\le r$. The scheme $\GL_r$
is its open subscheme obtained by localization by the determinant
polynomial $\det$. This is a group scheme (see \cite{Waterhouse}).
The group scheme $\SL_r=\Spec \Z[a_{i,j}]/(\det-1)$ is a closed
subscheme of $M_r.$ For every field $K$ the group schemes
$\GL_r(K)$ and $\SL_r(K)$ are obtained from $\GL_r$ and $\SL_r$ by
the base change. Then the groups $\GL_r(K)$ and $\SL_r(K)$ are the
groups of the $K-$rational geometric points of the corresponding
group schemes.

The multiplicative abelian group scheme $T_m$ acts on $M_r$ by
scalar multiplication. The scheme $\GL_r$ is invariant under this
action. This induces the action of the multiplicative group $K^*$
on the group $GL_r(K).$ The quotient of $\GL_r(K)$ by this action
is the group $\PGL_r(K).$

For every group word $w$ we consider the formal expression $\bar
w$ which is obtained from $w$ by replacing every letter $x\iv$ by
the symbol $\adj(x)$. Thus to every word $w$ in $k$ letters, we
can associate a polynomial map $\pi_w\colon M_r^k\to M_r$ which
takes every $k$-tuple of matrices $(A_1,...,A_k)$ to $\bar
w(A_1,...,A_k)$ where $\adj(A_i)$ is interpreted as the adjoint of
$A_i$. This map coincides with $w$ on $\SL_r$ since for the matrices
in $\SL_r$, the adjoint coincides with the inverse.

Similarly, for every endomorphism $\phi$ of the free group $F_k$,
we can extend the map $\phi_{\SL_r}\colon \SL_r^k\to \SL_r^k$ to a
self-map of $M_r^k$ which we shall denote by $\Phi$.

The map $\Phi$ is a self-map of the scheme $\Spec
\Z[a_{i,j}^{m}]$, $1\le i,j\le r, 1\le m\le k$. By base change it
induces a self-map of the scheme $\Spec K[a_{i,j}^{m}]$. $1\le
i,j\le r, 1\le m\le k$ for every field $K$. This map can be
restricted to the self-map of $\GL_r^k(K)$. The induced map on the
$K-$rational points is $\phi_{\GL_r(K)}$. It descends to the
self-map of $\PGL_r^k(K)$ which coincides with the map
$\phi_{\PGL_r(K)}$ defined above.

\medskip

Now we are ready to derive Theorems \ref{th} and \ref{thH} from
Theorems \ref{th0} and \ref{Hr}, respectively.

{\it Proof of Theorem \ref{th}.} Let $\phi$ be an injective
endomorphism of the free group $F_k=<x_1,...,x_k>$ and $1\neq w
\in F_k.$ Consider the self-map $\Phi$ of the scheme $M_2^k$ as
above. Denote $n=4k$. Similarly to Theorem \ref{th0}, we denote by
$V$ the Zariski closure of $\Phi^n(M_2^k).$ This is a scheme over
$\Spec{\Z}.$ Consider the map $\pi_w \colon V \to M_2$ as above.
We have the following lemma.

\begin{lemma} \label{sanov} In the above notations, $\pi_w(V)$
is not contained in the scheme of the scalar matrices.
\end{lemma}
\proof By the result of Sanov \cite{Sanov} there is an embedding
$\gamma \colon F_k \to \SL_2({\Z}).$ Obviously, $\gamma(F_k)$ does
not contain the matrix $-Id.$ Since $\phi $ is injective,
$\phi^n(w)\neq 1.$ Therefore
$\pi_w(\Phi^n(\gamma(x_1),...,\gamma(x_k))=\gamma(\phi^n(w))$ is
not a scalar matrix.
\endproof

Now we fix a big enough prime $p$ and make a base change from $\Z$
to the finite field $\F_p.$ Slightly abusing the notation, we will
from now on denote by $\Phi$ and $\pi_w$ the maps of the
corresponding schemes over $\F_p.$ And $V$ will also denote the
corresponding scheme over $\F_p.$ From Lemma \ref{sanov}, for big
enough $p$, $\pi_w{V}$ is not contained in the scheme of scalar
matrices. Consider the subscheme $Z_w$ of $V$ which is the union
of the $\pi_w$-pullback of scalar matrices and the subscheme of
$V$ consisting of $k$-tuples where one of the coordinates is
singular. We have that $Z_w$ is a proper subscheme of of $V$. By
Theorem \ref{th0} there exists a point $h=(a_1,...,a_k)\in
V\setminus Z_w$ such that $\Phi(h)=(a_1^Q,...,a_k^Q)$ for some
$Q=p^s$. Then the powers of $\Phi$ take the point $h$ to
$(a_1^{Q^l},...,a_k^{Q^l})$, $l\ge 1$ (we use the fact that, in
characteristic $p$, the Frobenius  commutes with every polynomial
map defined over $\F_p$). Therefore some power of $\Phi$ fixes
$h$. In addition $\pi_w(h)$ is not a scalar matrix and each $a_i$
is not a singular matrix because $h\not\in Z_w$. Taking the factor
over the torus action, we get a point $h'$ in $\PGL_2^k(\F_{p^i})$
that is fixed by some power of the map $\Phi$ and such that
$w(h')\ne 1$ in $\PGL_2(\F_{p^i})$. Thus the group
$\PGL_2(\F_{p^i})$ and the point $h'$ satisfy all three conditions
of Lemma \ref{lm2}. Since $w\in F_k$ was chosen arbitrarily, the
group $\HNN_\phi(F_k)$ is residually finite. This completes the
proof of Theorem \ref{th}.

{\it Proof of Theorem \ref{thH}.} Suppose $G\subseteq \SL_r(K).$
Here $K$ is some field and $G=<x_1,...,x_k \mid R>.$ Let $U_G$ be
the representation scheme of the group $G$ in $\SL_{r}$, i.e. the
reduced scheme of $k$-tuples of matrices from $\SL_{r}$ satisfying
the relations from $R$. This is a scheme over $\Spec K.$ Suppose
$\phi$ is an injective endomorphism of $G$ and $1\neq w \in G.$ We
choose a representation of $\phi$ and $w$ using the words on
$x_1,...,x_k$ and consider the maps $\Phi$ and $\pi_w.$  Since
$\phi$ is an endomorphism of $G$, the representation subscheme
$U_G$ is invariant under $\Phi$. Obviously, for big enough $m$ the
map $\Phi$ is dominant on the subscheme $V,$ which is the Zariski
closure of $\Phi^m(U_G)$. Note that $V$ may be reducible because
$U_G$ may be reducible. Since $\phi$ is injective, $\phi^m(w)\neq
1.$ Therefore $\pi_w(V)\neq \{Id\}$. By the usual specialization
argument (as in \cite{Mal}) there exists a finite field $\F_q$
such that for the corresponding schemes and maps over $\F_q$ the
same properties are satisfied. That is, $\Phi$ is dominant on $V,$
where $V$ is the Zariski closure of $\Phi^m(U_G)$,everything over
$\F_q$. In addition, $\pi_w(V)\neq \{Id\}$. Consider the subscheme
$Z_w$ of $V$ which is the $\pi_w$-pullback of the identity. We
have that $Z_w$ is a proper subscheme of $V$. We enlarge the
finite field to make all irreducible components of $V$ defined
over $\F_q.$  Some power of $\Phi$ maps each of these components
into itself, dominantly. We now apply Theorem \ref{Hr} to this
power of $\Phi,$ the scheme $V\subseteq \SL_r^k$ and its
subscheme $Z_w$. As in the proof of Theorem \ref{th}, we find a
point $h\in V(\F_{q^i})\setminus Z_w(\F_{q^i})$ that is fixed by
some power of $\Phi$. Since $h$ belongs to the representation
variety $U_G$, its coordinates satisfy all the relations from
$R$, so the condition (i) of Lemma \ref{lm2} is satisfied. Other
conditions of the lemma hold as before. Thus we can take the
group $\SL_r(\F_{q^i})$ as $H_w$, and $h$ as the point required
by Lemma \ref{lm2}. This completes the proof of Theorem \ref{thH}.

\section{Polynomial maps over finite fields}
In this section, we shall give a self-contained proof of Theorem
\ref{th0}.

Let $A^n$ be the affine space. Consider a map $\Phi:A^n
\rightarrow A^n$ given in coordinates by polynomials
$$f_1(x_1,...,x_n), f_2(x_1,...,x_n), ...f_n(x_1,...,x_n).$$
 The coordinate functions of the composition power $\Phi^k$ will be denoted by
$f_i^{(k)}(x_1,...,x_n),$ for $1\leq i\leq n.$ In what follows,
$\Phi$ will be defined over the finite field $\F_q$ of $q$
elements (this just means that all coefficients of $f_i$ belong to
$\F_q$). The number $Q$ will always mean some (big enough) power
of $q.$

We define by induction a chain of irreducible closed subvarieties
of $A^n$. Let $V_0=A^n,$ and for every $i\geq 1$ let $V_i$ be the
Zariski closure in $A^n$ of $\Phi(V_{i-1})$. Alternatively, $V_i$
is the Zariski closure of $\Phi^{i}(A^n).$

The varieties $V_i$ are irreducible (as polynomial images of an
irreducible variety) and $V_{i+1}\subseteq V_i$ for all $i$.
Because the dimension could only drop $n$ times, $V_{n}=
V_{n+1}=...$. We will denote this variety $V_n$ by $V$. Note that
$V=A^n$ if and only if $\Phi$ is a dominant map.

Suppose a point $a=(a_1,a_2,...,a_n) \in A^n$ is defined over the
algebraic closure $\overline{\F}_q$ of $\F_q.$ Recall that $a$ is called
quasi-fixed (with respect to $\Phi$) if there exists $Q=q^m$ such
that $f_i(a_1,a_2,...,a_n)=a_i^Q, i=1,...,n.$

In other words, the quasi-fixed points are those that are mapped
by $\Phi$ to their conjugates. They correspond to the closed
scheme points of $A^n,$ which are fixed by $\Phi$.

The following lemma is the first part of Theorem \ref{th0}.

\begin{lemma} \label{part1} All quasi-fixed points belong to the variety $V$.
\end{lemma}

\proof Since $\Phi$ is defined over $\F_q,$ all varieties $V_i,$
$i=1,2,...,n$ are defined over $\F_q$. For a point
$a=(a_1,...,a_n)\in A^n$ we denote by $a^Q$ the point
$\Fr_{q}^{m}(a)=(a_1^Q,...,a_n^Q).$ Then suppose $\Phi(a)=a^Q,$
for $Q=q^{m}.$ This implies that $a^Q\in V_1.$ Since the Frobenius
$\Fr_q$ commutes with $\Phi$, all varieties $V_i$ are invariant
with respect to $\Fr_q$. Therefore $a\in V_1.$ Hence
$a^Q=\Phi(a)\in V_2$ and $a\in V_2$. By induction, we get $a\in
V$.
\endproof

In the above notations, our main goal is to prove the following
(this is the second part of Theorem \ref{th0}).

\begin{theo} \label{mainag} Let $V$ be the Zariski closure of
$\Phi^n(A^n)$. Then quasi-fixed points of $\Phi$ are Zariski dense
in $V$. In other words, suppose $W\subset V$ is a proper Zariski
closed subvariety. Then for some $Q$ there is a point
$(a_1,...,a_n)\in V\setminus W$ such that $f_i(a_1,...a_n)=a_i^Q,
i=1,...,n.$
\end{theo}

We denote by $I_Q$ the ideal in $\F_Q[x_1,...,x_n]$ generated by
the polynomials $f_i(x_1,...,x_n)-x_i^Q,$ for $i=1,2,...,n.$

\begin{lemma} \label{lm5} For a big enough $Q$ the ideal $I_Q$ has finite length.
\end{lemma}
\proof We compactify $A_n$ to the projective space $P^n$ in the
usual way. We also projectivize the polynomials
$f_i-x_i^Q$. If there is a curve in $P^n$ on which all of these
projective polynomials vanish, then it must have some points on
the infinite hyperplane of $P^n$. But this is impossible if $Q$ is
bigger than the degrees of $f_i.$. Thus the scheme of common
zeroes is zero-dimensional, which implies the result.
\endproof

\begin{lemma} \label{lm6} For all $1\leq i\leq n$ and $j\geq 1$
$$f_i^{(j)}(x_1,...,x_n)-x^{Q^j}\in I_Q.$$
\end{lemma}
\proof We use induction on $j$. For $j=1$ the statement is
obvious. Suppose it is true for some $j\ge 1.$ Then
$$f_i^{(j+1)}(x_1,...,x_n)=f_i(f_1^{(j)},...,f_n^{(j)}) \equiv f_i(x_1^{Q^j},...,x_n^{Q^j})=$$
$$=f_i(x_1,...,x_n)^{Q^j}\equiv x^{Q^{j+1}} (\mod I_Q)$$
\endproof

The next lemma is the crucial step in the proof.
\begin{lemma}\label{lm7}
There exists a number $k$ such that for every quasi-fixed point
$(a_1,...,a_n)$ with big enough $Q$ and for every $1\leq i\leq n$
the polynomial
$$(f_i^{(n)}(x_1,...,x_n)- f_i^{(n)}(a_1,...,a_n))^k$$
is contained in the localization of $I_Q$ at $(a_1,...,a_n).$
\end{lemma}
\proof Let us fix $i$ from $1$ to $n$. The polynomials $x_i, f_i,
f_i^{(2)},..., f_i^{(n)}$ are algebraically dependent over $\F_q.$
This means that \begin{equation}\label{eq67}\sum \limits_{s}
a_s\cdot (x_i)^{\alpha_{0,s}} \cdot (f_i)^{\alpha_{1,s}}\cdot
...\cdot (f_i^{(n)})^{\alpha_{n,s}} =0 \end{equation} with some
non-zero $a_s\in \F_q.$ By Lemma \ref{lm6} the polynomial in the
left hand side of (\ref{eq67}) is congruent modulo $I_Q$ to
$$P_Q(x_i)=\sum \limits_{s} a_s\cdot x_i^{\alpha_s},$$
where $\alpha_s=\sum \limits_{j=0}^{n}\alpha_{j,s}Q^j.$ For big
enough $Q$, the polynomial $P_Q$ is non-zero. For any $(a_1,...a_n)$
we rewrite $P_Q(x_i)$ as $\sum \limits_{t} b_t\cdot
(x_i-a_i)^{\beta_t}.$

So in the local ring of $(a_1,...,a_n)$, the polynomial $P_Q(x_i)$
is equal to $$(x_i-a_i)^{\beta} \cdot u,$$ where $u$ is invertible
and $\beta \leq \max \beta_t.$ Clearly, $\max \beta_t$ is bounded
by $kQ^n$ for some $k$ that does not depend on $Q, a_1,...,a_n$.
Denote by $I_Q^{(a_1,...,a_n)}$ the localization of $I_Q$ in the
local ring of $(a_1,...,a_n).$ Then by (\ref{eq67})
$(x_i-a_i)^{kQ^n}\equiv 0 \quad (\mod I_Q^{(a_1,...,a_n)}).$ Now
we note that
$$f_i^{(n)} (x_1,...,x_n) - f_i^{(n)} (a_1,...,a_n) =$$
$$=f_i^{(n)} (x_1,...,x_n) -a_i^{Q^n} \equiv x_i^{Q^n}-a_i^{Q^n}=(x_i-a_i)^{Q^n} (\mod I_Q^{(a_1,...,a_n)}). $$
\endproof

Let us fix some polynomial $D$ with the coefficients in a finite
extension of $\F_q$ such that it vanishes on $W$ but not on $V$.
By base change we will assume that all coefficients of $D$ are in
$\F_q.$

\begin{lemma} \label{P} There exists a positive integer $K$ such that for all quasi-fixed
points $(a_1,...,a_n)\in W$ with big enough $Q$ we get
$$(D(f_1^{(n)}(x_1,...,x_n),...,f_n^{(n)}(x_1,...,x_n)))^K\equiv 0 (\mod I_Q^{(a_1,...,a_n)})$$
\end{lemma}

\proof For every $(a_1,...,a_n)\in W$ we can rewrite
$D(x_1,...,x_n)$ as a polynomial in $x_i-a_i^{Q^n}.$ This
polynomial has no free term because $D$ vanishes on $W$ and
$(a_1,...,a_n)\in W$ by the assumption. The number of non-zero
terms of $D$ is bounded independently of $a_i$ and $Q$ by some
number $N$. Then by the binomial formula and Lemma \ref{lm7} we
can take $K=N(k-1)+1$ where $k$ is the constant from Lemma
\ref{lm7}.
\endproof

The polynomial
$P=(D(f_1^{(n)}(x_1,...,x_n),...,f_n^{(n)}(x_1,...,x_n)))^K$ is
non-zero because $D$ does not vanish on the whole $U$ and the map
$\Phi$ is dominant on $U$. We now complete the proof of
Theorem \ref{mainag}.

In fact we give two proofs. The first one uses the Bezout theorem,
while the second one is elementary and self-contained.

{\bf Proof \# 1.} We denote by $Z$ the subscheme of $A^n$
that corresponds to $P$. Note that $Z$ does not depend on $Q.$ Now
for every $Q$ consider the $\F_q-$linear subspace of polynomials
spanned by $f_i-x_i^Q,$ $1\leq i \leq n.$ By Lemma \ref{lm5} its
base locus is zero-dimensional, i.e. these polynomials do not
vanish simultaneously at any curve. The scheme $Z$ is of pure
dimension $(n-1)$. A general element $\tau_1$ of the above linear
subspace does not vanish at any of the irreducible components of
$Z,$ or their positive-dimensional intersections. So its scheme of
zeroes intersects $Z$ properly,  the intersection $Z_1$ has pure
dimension $(n-2).$ Then we choose $\tau_2$ that intersects $Z_1$
properly to get $Z_2,$ and so on. After choosing $(n-1)$ elements
$\tau_1, \tau_2,...,\tau_{n-1}$ we get an ideal $I_Q'\la D^k,
\tau_1, \tau_2,...,\tau_{n-1}\ra$ of finite length. After
localization at any $(a_1,...a_n)\in W$ this ideal is contained in
$I_Q.$ By Bezout theorem (cf., e.g. \cite{Fulton}) )  the length
of $I_Q'$ is equal to $const\cdot Q^{n-1}.$ But the length of
$I_Q$ is equal to $Q^n,$ which is bigger for big enough $Q.$ This
implies the existence of quasi-fixed points in $V\setminus W.$
\endproof

{\bf Proof \# 2.} Fix $Q=q^i$ such that it is bigger than the
degrees of $f_i$ and $P$.  By Lemma \ref{part1} all points with $\Phi(x)=x^Q$ belong to $V$. If they all actually
belong to $W$, then $P$ lies in the localizations of $I_Q$ with
respect to all maximal ideals of $\F[x_1,..x_n]$. Therefore $P\in
I_Q$ (otherwise consider a maximal ideal containing $I_Q:F$). This
means that there exist polynomials $u_1,...u_n$ such that
\begin{equation}\label{eq90}
P=\sum \limits_{i=1}^{n} u_i \cdot (f_i-x_i^Q)
\end{equation}
The right hand side of (\ref{eq90}) can be modified as follows.
For every $i<j$ and any polynomial $A$, we can add $A(f_j-x_j^Q)$
to $u_i$ and subtract $A(f_i-x_i^Q)$ from $u_j$. This
transformation can be used to make (for every $i<j$) the degree of
$x_i$ in every monomial in $u_j$ is smaller than $Q$. Now consider
the monomial of the highest total power in (\ref{eq90}). Since $Q$
is bigger than the degree of $P$, that monomial has the form $\bar
u_jx_j^Q$ for some $j$ where $\bar u_j$ is the leading monomial in
$u_j$. That monomial does not occur in the left hand side of
(\ref{eq90}). Therefore it must coincide with the leading monomial
$\bar u_ix_i^Q$ for some $i\ne j$. But then $u_i$ must be
divisible by $x_j^Q$ and $u_j$ must be divisible by $u_i^Q$, and
we get a contradiction in each of the cases $i<j$ or $i>j$.
\endproof

\section{Extendable endomorphisms of linear\\ groups, and some open problems}

Recall that a profinite group is, by definition, a projective
limit of finite groups.

\begin{df}{\rm
Let $G$ be a residually finite group, $\phi$ be an endomorphism of
$G$. We say that $\phi$ is {\em extendable} if there exists a
profinite group $\bar G$  containing $G$ as a dense subgroup, and
a (continuous) automorphism $\bar\phi$ of $\bar G$ such that
$\phi$ is the restriction of $\bar\phi$ on $G$. }\end{df}

Notice that even if $\phi$ is injective and continuous in a
profinite topology of $G$, its (unique) extension to the
corresponding completion of $G$ may not be injective. Injective
endomorphisms of free groups that have injective extensions in
$p$-adic (resp. pro-solvable, and many other profinite) topologies
of a free group are completely described in \cite{CSW}.

There is  a close connection between extendable endomorphisms and
residually finite HNN extensions.

\begin{theo} \label{thext} An injective endomorphism $\phi$ of a residually
finite group $G$ is extendable if and only if $\HNN_\phi(G)$ is
residually finite.
\end{theo}

\proof Suppose that $P=\HNN_\phi(G)$ is residually finite. Let
$\Psi$ be the set of all homomorphisms of $P$ onto finite groups,
$\Psi'$ be the set of all restrictions of homomorphisms from
$\Psi$ to $G$. Let $\ttt$ be the smallest profinite topology on
$G$ for which all the homomorphisms from $\Psi'$ are continuous.
The base of neighborhoods of $1$ for $\ttt$ is formed by the
kernels of all the homomorphisms from $\Psi'$.

It is easy to see that for every $\psi\in\Psi'$, the homomorphism
$\phi\psi$ is also in $\Psi'$. Therefore the endomorphism $\phi$
is continuous in the topology $\ttt$. Let $\bar G$ be the
profinite completion of $G$ with respect to $\ttt$, and let
$\bar\phi$ be the (unique) continuous extension of $\phi$ onto
$\bar G$. Let us prove that $\bar\phi$ is an automorphism of $\bar
G$.

Suppose that $\bar\phi$ is not injective. This means that there is
a sequence of elements $w_i$, $i\ge 1$, in $G$ such that $w_i$ do
not converge to $1$ in $\bar G$ but $\phi(w_i)$ converge to $1$.
The latter means that there exists a sequence of subgroups
$N_i=\Ker(\psi_i)\subset P$, $\psi_i\in \Psi$, such that $\cap
N_i=\{1\}$, $\phi(w_i)\in N_i$, $i\ge 1$.

Notice that by definition of $P=\HNN_\phi(G)$, $\phi(w_i)N_i$ is a
conjugate of $w_iN_i$ in $P/N_i$ (the conjugating element is
$tN_i$ where $t$ is the free letter of the HNN extension). Thus we
can conclude that $w_i\in N_i$, $i\ge 1$. Hence $w_i\to 1$ in
$\ttt$, a contradiction. Therefore $\bar\phi$ is injective.

Let us prove that $\bar\phi$ is surjective. Consider a Cauchy
sequence $w=\{w_i, i\ge 1\}$ in $G$, that is suppose there exist
$N_i=\Ker (\psi_i)$, $\psi_i\in \Psi$, $i\ge 1$, such that $\cap
N_i=\{1\}$ and $w_i\iv w_j\in N_i$ for every $j>i$.

For every $x\in G$ we have $\phi(x)N_i=txt\iv N_i$, and $P/N_i$ is
finite. So $\phi$ induces an automorphism in $G/(N_i\cap G)$. Hence for every $i\ge 1$, we can find an element $u_i$ in $G$ such that
$\phi(u_i)N_i=w_iN_i$. Moreover since $w_i\iv w_j\in N_i$ for all
$j>i$, $u_i\iv u_j\in N_i$ as well. Therefore $\{u_i, i\ge i\}$ is
a Cauchy sequence and $\bar\phi(u)=w$. Thus $\bar\phi$ is an
automorphism of $\bar G$. Notice that since $\bar G$ is compact,
$\phi\iv$ is also continuous.

Suppose now that $\phi$ can be extended to a continuous
automorphism $\bar\phi$ of a profinite group $\bar G\ge G$. Let
$w\ne 1\in G$. Notice that for every $w\in G$,
$\bar\phi(w)=\phi(w)$. Therefore there exists a homomorphism
$\theta$ from $P$ to the semidirect product $\bar G\rtimes
\la\bar\phi\ra$ which is identity on $G$ and sends $t$ to
$\bar\phi$. This homomorphism is clearly injective: it is easy to
check that no non-trivial element $t^k w t^l$ can lie in the
kernel of $\theta$. It remains to prove that $\bar G\rtimes
\la\bar\phi\ra$ is residually finite. But that can be done exactly
as in the case of split extensions of finitely generated
residually finite groups \cite{Mal1}. Indeed, since $\bar G$ is
finitely generated as a profinite group, it has only finitely many
open subgroups of any given (finite) index, and the automorphism
$\bar\phi$ permutes these subgroups. Hence $\bar\phi$ leaves
invariant their intersection which also is of finite index.
Therefore $\bar G\rtimes \la\phi\ra$ is residually
finite-by-cyclic, so $G\rtimes \la\phi\ra$ is residually
finite.\endproof

Theorems \ref{thH} and \ref{thext} immediately imply:

\begin{cy} Every injective endomorphism of a finitely generated
linear group is extendable.
\end{cy}

Finally let us mention two open problems.

\begin{prob} {\rm Let $\phi$ and $\psi$ be two injective endomorphisms
of the free group $F_k=\la x_1,...,x_k\ra$. Consider the
corresponding HNN extension of $F_k$ with two free letters $t,u$:
$$\HNN_{\phi,\psi}(F_k)=\la x_1,...,x_k, t, u\mid
tx_it\iv=\phi(x_i), ux_iu\iv=\psi(x_i), 1\le i\le k\ra.$$ Is
$\HNN_{\phi,\psi}(F_k)$ always residually finite? }\end{prob}

We believe that the answer is negative in a very strong sense: the
groups $\HNN_{\phi,\psi}$ should be generically non-residually
finite. Since many of these groups are hyperbolic, this may
provide a way to construct hyperbolic non-residually finite
groups.

The next question is natural to ask for any residually finite
groups.

\begin{prob}{\rm Are mapping tori of free group endomorphisms linear?
}\end{prob}

Notice that although we use linear groups in our proof of Theorems
\ref{th} and \ref{thH}, our proof does not prove linearity of the
mapping tori. It is easy to extract from our proof that the
mapping torus a linear group endomorphism is embeddable into the
wreath product of a linear group and the infinite cyclic group.
Notice that this wreath product is not even residually finite.

\begin{minipage}[t]{2.5 in}
\noindent Alexander Borisov\\Department of Mathematics\\
Penn State University \\ borisov@math.psu.edu\\
\end{minipage}
\begin{minipage}[t]{2.5 in}
\noindent Mark V. Sapir\\
Department of Mathematics\\
Vanderbilt University\\
m.sapir@vanderbilt.edu\\
\end{minipage}

\end{document}